\documentclass[12pt,twoside]{amsart}
\usepackage{amssymb}
\usepackage{amscd}

\title[Log canonical ring in dimension four]
{Finite generation of the log canonical 
ring in dimension four}
\author{Osamu Fujino}
\dedicatory{Dedicated to the memory of Professor Masayoshi Nagata}
\date{2009/9/23, version 5.01}
\address{Department of Mathematics, Faculty of Science, 
Kyoto University, Kyoto 606-8502 Japan}
\email{fujino@math.kyoto-u.ac.jp}

\newcommand{\Supp}[0]{{\operatorname{Supp}}}

\newcommand{\Exc}[0]{{\operatorname{Exc}}}

\newcommand{\Alb}[0]{{\operatorname{Alb}}}

\newtheorem{thm}{Theorem}[section]
\newtheorem{lem}[thm]{Lemma}
\newtheorem{cor}[thm]{Corollary}
\newtheorem{prop}[thm]{Proposition}

\newtheorem{cla}{Claim}
\newtheorem{conj}[thm]{Conjecture}

\theoremstyle{definition}
\newtheorem{ex}[thm]{Example}
\newtheorem{defn}[thm]{Definition}
\newtheorem{rem}[thm]{Remark}
\newtheorem*{ack}{Acknowledgments}       

\newtheorem{saya}{}
\newtheorem{prob}[thm]{Problem}
\begin{document}
\bibliographystyle{amsalpha+}

\begin{abstract}
We treat two different topics on the log minimal model program, 
especially for four-dimensional log canonical pairs. 

(a) Finite generation of the log canonical ring in dimension four. 

(b) Abundance theorem for irregular fourfolds. 

We obtain (a) as a direct consequence of 
the existence of four-dimensional log 
minimal models by using Fukuda's 
theorem on the four-dimensional log abundance conjecture. 
We can prove (b) only by using traditional arguments. 
More precisely, we prove the abundance conjecture for irregular 
$(n+1)$-folds on the assumption that 
the minimal model conjecture and 
the abundance conjecture hold in dimension $\leq n$. 
\end{abstract}

\maketitle
\tableofcontents

\section{Introduction}\label{sec1}  

In this paper, we treat two different topics on 
the log minimal model program, especially for 
four-dimensional log canonical pairs. We will freely use 
the results on the three-dimensional 
log minimal model program (cf.~\cite{fafa}, 
\cite{kemm}, etc.). 
Sorry, we do not always refer to the original 
papers since the results are scattered in various places. 

\begin{saya}
[{\textbf{Finite generation of the log canonical ring in dimension four}}] 

The following theorem is the main result of Section \ref{sec3} 
(cf.~\cite[Section 3.1]{fujino-kmm}). 

\begin{thm}[Finite generation of the 
log canonical 
ring in dimension four]\label{main}
Let $\pi:X\to Z$ be a proper surjective morphism 
from a smooth fourfold $X$. 
Let $B$ be a boundary $\mathbb Q$-divisor 
on $X$ such that $\Supp B$ is a simple normal crossing 
divisor on $X$. 
Then the relative log canonical ring 
$$
R(X/Z, K_X+B)=
\bigoplus _{m\geq 0}\pi_*\mathcal O_X(\llcorner m(K_X+B)\lrcorner)
$$ 
is a finitely generated $\mathcal O_Z$-algebra. 
\end{thm}

It is easy to see that 
Theorem \ref{main} is equivalent to 
Theorem \ref{main2}. 

\begin{thm}\label{main2}
Let $\pi:X\to Z$ be a proper surjective morphism 
from a four-dimensional log canonical pair $(X, B)$ such 
that $B$ is an effective $\mathbb Q$-divisor. 
Then the relative log canonical ring 
$$
R(X/Z, K_X+B)=
\bigoplus _{m\geq 0}\pi_*\mathcal O_X(\llcorner m(K_X+B)\lrcorner)
$$ 
is a finitely generated $\mathcal O_Z$-algebra. 
\end{thm}

In Section \ref{sec3}, we give a proof of 
Theorem \ref{main} by using the existence theorem of 
four-dimensional log minimal models (cf.~\cite{birkar} and \cite{shokurov2}) and 
Fukuda's result on the log abundance conjecture for 
fourfolds (cf.~\cite{fukuda}). A key point of 
Fukuda's result is the abundance theorem 
for semi log canonical threefolds in \cite{fujino1}. 
\end{saya}

\begin{saya}[{\textbf{Abundance theorem for irregular fourfolds}}] 
In Section \ref{sec4}, we prove the 
abundance theorem for irregular $(n+1)$-folds 
on the assumption that 
the minimal model conjecture and 
the abundance conjecture hold in dimension $\leq n$ 
(see Theorem \ref{thm-abun}). 
By this result, we know that 
the abundance conjecture for irregular varieties is the problem 
for lower dimensional varieties. 
Since the minimal model conjecture and the abundance conjecture 
hold in dimension $\leq 3$, 
we obtain the next theorem (see Corollary \ref{cor45}). 

\begin{thm}[Abundance theorem for irregular fourfolds]
Let $X$ be a normal complete fourfold with only canonical 
singularities. 
Assume that $K_X$ is nef and the irregularity $q(X)$ is not zero. 
Then $K_X$ is semi-ample. 
\end{thm}
We also prove that 
there exists a good minimal model for 
any smooth projective irregular fourfold (see Theorem \ref{thm-good}).  

\begin{thm}[Good minimal models of irregular fourfolds]
Let $X$ be a smooth 
projective irregular fourfold. 
If $X$ is not uni-ruled, then 
there exists a normal 
projective variety $X'$ such that 
$X'$ has only $\mathbb Q$-factorial 
terminal singularities, $X'$ is birationally equivalent to 
$X$, and $K_{X'}$ is semi-ample.  
\end{thm}
\end{saya}
We note that Sections \ref{sec3} and \ref{sec4} can be read independently. 

\begin{ack}
The author was partially supported 
by the Grant-in-Aid for Young Scientists (A) $\sharp$20684001 from 
JSPS. 
He was also supported by the Inamori Foundation.  
\end{ack}

We will work over $\mathbb C$, the complex number field, 
throughout this paper. 
We will freely use the notation in \cite{kmm} and \cite{km}. 
Note that we do not use $\mathbb R$-divisors. 

\section{Preliminaries}\label{sec2}
In this section, we collect basic definitions. 

\begin{defn}[Divisors, $\mathbb Q$-divisors] 
Let $X$ be a normal variety. 
For a $\mathbb Q$-Weil divisor 
$D=\sum _{j=1}^r d_j D_j$ on $X$ such that 
$D_j$ is a prime divisor for every 
$j$ and $D_i\ne D_j$ for $i\ne j$, we define 
the {\em{round-down}} $\llcorner D\lrcorner 
=\sum _{j=1}^{r} \llcorner d_j \lrcorner D_j$, 
where for every rational number $x$, 
$\llcorner x\lrcorner$ is the integer defined by 
$x-1<\llcorner x\lrcorner \leq x$. 

We call $D$ a {\em{boundary}} 
$\mathbb Q$-divisor if 
$0\leq d_j\leq 1$ 
for every $j$. 

We note that $\sim _{\mathbb Q}$ 
denotes the $\mathbb Q$-linear equivalence 
of $\mathbb Q$-divisors. 

We call $X$ {\em{$\mathbb Q$-factorial}} if and only if 
every Weil divisor on $X$ is $\mathbb Q$-Cartier. 
\end{defn}

\begin{defn}[Exceptional locus]
For a proper birational morphism $f:X\to Y$, 
the {\em{exceptional locus}} $\Exc (f)\subset X$ is the locus where 
$f$ is not an isomorphism. 
\end{defn}

Let us quickly recall the definitions of singularities of pairs. 

\begin{defn}[Singularities of pairs] 
Let $X$ be a normal variety and $B$ an effective $\mathbb Q$-divisor 
on $X$ such that $K_X+B$ is $\mathbb Q$-Cartier. 
Let $f:Y\to X$ be a resolution such that 
$\Exc (f)\cup f^{-1}_*B$ has a simple normal crossing 
support, where $f^{-1}_*B$ is the strict transform of $B$ on $Y$. 
We write $$K_Y=f^*(K_X+B)+\sum _i a_i E_i$$ and 
$a(E_i, X, B)=a_i$. 
We say that $(X, B)$ is {\em{lc}} (resp.~{\em{klt}}) if and only if 
$a_i\geq -1$ (resp.~$a_i>-1$) for every $i$. 
We note that lc (resp.~klt) is an abbreviation of 
{\em{log canonical}} (resp.~{\em{Kawamata log terminal}}). 
We also note that the {\em{discrepancy}} $a(E, X, B)\in \mathbb Q$ can be 
defined for every prime divisor $E$ {\em{over}} $X$. 

In the above notation, if $B=0$ and $a_i>0$ (resp.~$a_i\geq 0$) 
for every $i$, then we 
say that $X$ has only {\em{terminal}} 
(resp.~{\em{canonical}}) singularities. 
\end{defn}

\begin{defn}[Divisorial log terminal pair] 
Let $X$ be a normal variety and 
$B$ a boundary $\mathbb Q$-divisor 
such that 
$K_X+B$ is $\mathbb Q$-Cartier. 
If there exists a resolution 
$f:Y\to X$ 
such that 
\begin{itemize}
\item[(i)] both $\Exc (f)$ and $\Exc (f)\cup \Supp (f^{-1}_*B)$ 
are simple normal crossing divisors on $Y$, and 
\item[(ii)] $a(E, X, B)>-1$ for every 
exceptional divisor $E\subset Y$, 
\end{itemize}
then $(X, B)$ is called {\em{divisorial log terminal}} 
({\em{dlt}}, for short). 
\end{defn}

For the details of singularities of pairs, 
see, for example, \cite{km} and \cite{what}. 

\begin{defn}[Center, lc center] 
Let $E$ be a prime divisor over $X$. The closure 
of the image of $E$ on $X$ is denoted by 
$c_X(E)$ and 
called the {\em{center}} of $E$ on $X$. 

Let $(X, B)$ be an lc pair. 
If $a(E, X, B)=-1$, $c_X(E)$ is called 
an {\em{lc center}} of $(X, B)$. 
\end{defn} 

The following definitions are now classical. 

\begin{defn}[Iitaka's $D$-dimension and numerical 
$D$-dimension]
Let $X$ be a normal complete  
variety and $D$ a $\mathbb Q$-Cartier $\mathbb Q$-divisor. 
Assume that $mD$ is Cartier for a positive 
integer $m$. 
Let 
$$\Phi_{|tmD|}: X\dashrightarrow \mathbb P^{\dim |tmD|}
$$ 
be rational mappings given by 
linear systems $|tmD|$ for positive integers $t$. 
We define {\em{Iitaka's $D$-dimension}} 
\begin{eqnarray*}
\kappa (X, D)=\left\{
\begin{array}{ll} 
\underset{t>0}{\max} \dim \Phi _{|tmD|}(X), & {\text{if}} \ \ 
|tmD|\ne \emptyset 
\ \ {\text{for some}}\ \  t, 
\\ 
-\infty, & {\text{otherwise}}.
\end{array}\right. 
\end{eqnarray*}
In case $D$ is nef, 
we can also define the {\em{numerical 
$D$-dimension}} 
$$
\nu(X, D)=\max \{ \, e \,|\, D^e\not \equiv 0 \}, 
$$ 
where $\equiv$ denotes {\em{numerical equivalence}}. 
We note that $\nu(X, D)\geq \kappa (X, D)$ always holds. 
\end{defn}

\begin{defn}[Nef and abundant divisors]\label{defn1}
Let $X$ be a normal complete variety and $D$ a $\mathbb Q$-Cartier 
$\mathbb Q$-divisor 
on $X$. 
Assume that $D$ is nef. 
The nef $\mathbb Q$-divisor $D$ is said to 
be {\em{abundant}} if the equality $\kappa (X, D)=
\nu(X, D)$ holds. 
Let $\pi:X\to Z$ be a proper surjective morphism of normal 
varieties and $D$ a $\pi$-nef 
$\mathbb Q$-divisor on $X$. Then $D$ is 
said to be {\em{$\pi$-abundant}} 
if $D_\eta$ is abundant, where 
$D_\eta=D|_{X_\eta}$ and 
$X_\eta$ is the generic fiber of $\pi$.   
\end{defn}

\begin{defn}[Irregularity] 
Let $X$ be a normal complete variety with only rational singularities. 
We put 
$$q(X)=h^1(X, \mathcal O_X)=\dim H^1(X, \mathcal O_X)<\infty$$ 
and call it the {\em{irregularity}} of $X$. 

Let $X$ be as above. 
If $q(X)\ne 0$, then 
we call $X$ {\em{irregular}}.

If $X'$ is a normal complete variety with only rational 
singularities such that 
$X'$ is birationally equivalent to 
$X$, then it is easy to 
see that $q(X)=q(X')$. 
 
\end{defn}

\section{Log canonical ring}\label{sec3}

In this section, we prove the following theorem:~Theorem \ref{main}. 

\begin{thm}[Finite generation of the four-dimensional 
log canonical 
ring]
Let $\pi:X\to Z$ be a proper surjective morphism 
from a smooth fourfold $X$. 
Let $B$ be a boundary $\mathbb Q$-divisor 
on $X$ such that $\Supp B$ is a simple normal crossing 
divisor on $X$. 
Then the relative log canonical ring 
$$
R(X/Z, K_X+B)=
\bigoplus _{m\geq 0}\pi_*\mathcal O_X(\llcorner m(K_X+B)\lrcorner)
$$ 
is a finitely generated $\mathcal O_Z$-algebra. 
\end{thm}

The next proposition is well known and a slight generalization of \cite[Theorem 
7.3]{kawamata2}. 

\begin{prop}\label{prop-a}
Let $(X, B)$ be a proper log canonical fourfold such that 
$K_X+B$ is nef and $\kappa (X, K_X+B)>0$. 
Then $K_X+B$ is abundant, that is, 
$\kappa (X, K_X+B)=\nu (X, K_X+B)$. 
\end{prop}
\begin{proof}
See, for example, \cite[Proposition 3.1]{fukuda}. 
We note that we need the three-dimensional 
log minimal model 
program and 
log abundance theorem here (see ~\cite{fafa}, \cite{kemm}, and \cite{kemm2}). 
\end{proof}

Let us recall Fukuda's result in \cite{fukuda}. 
We will generalize this in Theorem \ref{final}. 

\begin{thm}[{cf.~\cite[Theorem 1.5]{fukuda}}]\label{thm3}
Let $(X, B)$ be a proper dlt fourfold. 
Assume that $K_X+B$ is nef and that 
$\kappa (X, K_X+B)>0$. 
Then $K_X+B$ is semi-ample. 
\end{thm}
\begin{proof}
By Proposition \ref{prop-a}, $\kappa (X, K_X+B)=\nu (X, K_X+B)$. 
We put $S=\llcorner B\lrcorner$ and 
$K_S+B_S=(K_X+B)|_S$. Then the pair 
$(S, B_S)$ is semi divisorial 
log terminal and $K_S+B_S$ is semi-ample 
by \cite[Theorem 0.1]{fujino1}. Finally, by \cite[Corollary 6.7]{fujino2}, 
we obtain that $K_X+B$ is semi-ample. 
\end{proof}

\begin{rem}
The proof of \cite[Proposition 3.3]{fukuda} depends on 
\cite[Theorem 5.1]{kawamata2}. It 
requires \cite[Theorem 4.3]{kawamata2} whose 
proof contains a nontrivial gap. 
See \cite[Remark 3.10.3]{what} and 
\cite{fujino-kawamata}. 
So, we adopted \cite[Corollary 6.7]
{fujino2} in the proof of Theorem \ref{thm3}. 
\end{rem}

In this section, 
we adopt Birkar's definition of the log minimal model 
(see \cite[Definition 2.4]{birkar}), which is 
slightly different from \cite[Definition 3.50]{km}. 
See Remark \ref{rere} and Example \ref{rem24} below. 

\begin{defn}[{cf.~\cite[Definition 2.4]{birkar}}]
\label{lmm-def}  
Let $(X, B)$ be a log canonical 
pair over $Z$. 
A {\em{log minimal model}} $(Y/Z, B_Y+E)$ of 
$(X/Z, B)$ consists of a birational 
map $\phi: X\dashrightarrow Y/ Z$, $B_Y=\phi_*B$, and 
$E=\sum _j E_j$, where 
$E_j$ is a prime divisor on $Y$ and $\phi^{-1}$-exceptional for every $j$, and satisfies 
the following conditions: 
\begin{itemize}
\item[(1)] $Y$ is $\mathbb Q$-factorial 
and $(Y, B_Y+E)$ is 
dlt, 
\item[(2)] $K_Y+B_Y+E$ is nef over $Z$, and 
\item[(3)] for every prime divisor $D$ on $X$ which 
is 
exceptional over $Y$, 
we have 
$$
a(D, X, B)<a(D, Y, B_Y+E), 
$$
where $a(D, X, B)$ (resp.~$a(D, Y, B_Y+E)$) 
denotes the discrepancy of $D$ with respect to 
$(X, B)$ (resp.~$(Y, B_Y+E)$). 
\end{itemize}
\end{defn}

\begin{rem}\label{rere}
In \cite[Definition 3.50]{km}, 
it is required 
that $\phi^{-1}$ has no exceptional divisors. 
\end{rem}

\begin{ex}\label{rem24}  
Let $X=\mathbb P^2$ and $D_X$ the complement 
of the big torus. Then 
$K_X+D_X$ is dlt and $K_X+D_X\sim 0$. 
Let $Y=\mathbb P_{\mathbb P^1}(\mathcal O_{\mathbb P^1}
\oplus \mathcal O_{\mathbb P^1}(-1))$ and 
$D_Y$ the complement of the 
big torus. 
Then $(Y, D_Y)$ is a log minimal 
model 
of $(X, D_X)$ in the sense of Definition 
\ref{lmm-def}. 
Of course, $K_Y+D_Y$ is dlt and $K_Y+D_Y\sim 0$. 
On the 
other hand, $(Y, D_Y)$ is not a log minimal 
model of $(X, D_X)$ 
in the sense of \cite[Definition 3.50]{km}. 
\end{ex}

We prepare the following two easy lemmas. 

\begin{lem}\label{iso}
We use the notation in {\em{Definition \ref{lmm-def}}}. Then 
we have 
$$
a(\nu, X, B)\leq a(\nu, Y, B_Y+E)
$$ 
for every divisor $\nu$ over $X$. 
Thus, we obtain 
$$
R(X/Z, K_X+B)\simeq R(Y/Z, K_Y+B_Y+E). 
$$
\end{lem}
\begin{proof}
It is an easy consequence of the negativity lemma. 
See, for example, \cite[Proposition 3.51 and Theorem 3.52]{km}. 
\end{proof}

\begin{lem}\label{leng} 
Let $\pi:X\to Z$ be a projective surjective 
morphism between projective varieties. Assume that 
$(X, B)$ is log canonical 
and $H$ is an ample Cartier divisor on $Z$. 
Let $R$ be a $(K_X+B)$-negative 
extremal ray of $\overline {NE}(X)$ such that 
$$
R\cdot (K_X+B+(2\dim X+1)\pi^*H)<0. 
$$ 
Then $R$ is a $(K_X+B)$-negative 
extremal ray of $$\overline {NE}(X/Z)=\{ z \in \overline {NE}(X)\, 
| \, z\cdot \pi^*H=0\}\subset \overline {NE}(X). 
$$
In particular, if $K_X+B$ is $\pi$-nef, 
then $K_X+B+(2\dim X+1)\pi^*H$ 
is nef. 
\end{lem}
\begin{proof}
If $(X, B)$ is klt, then it is obvious by Kawamata's bound of the 
length of extremal rays (see \cite{kawamata3}). 
When $(X, B)$ is lc, it is sufficient to use \cite[Subsection 3.1.3]{fujino-kmm} or \cite[Section 18]{ff}. 
\end{proof}

Let us start the proof of 
Theorem \ref{main}. 

\begin{proof}[{Proof of {\em{Theorem \ref{main}}}}] 
We can assume that the fiber of $\pi$ is connected. 
First, if $\kappa (X_\eta, K_{X_{\eta}}+B_\eta)=-\infty$, 
where $\eta$ is the generic point of $Z$, then the 
statement is trivial. We note that 
the statement is obvious when $Z$ is a point and 
$\kappa (X, K_X+B)=0$. 
So, we can assume that $\kappa (X_\eta, K_{X_\eta}+ 
B_\eta)\geq 0$ and that 
$\kappa (X, K_X+B)\geq 1$ when $Z$ is a point. 
Since the problem is local, we can assume that $Z$ is affine. 
By compactifying $Z$ and $X$ and taking a resolution of $X$, we can assume 
that $X$ and $Z$ are projective and that 
$\Supp B$ is a simple normal crossing 
divisor. 
By the assumption, we can find an 
effective $\mathbb Q$-divisor $M$ on $X$ such that 
$K_X+B\sim _{\mathbb Q,\pi}M$, that is, 
there exists a $\mathbb Q$-divisor $N$ on $Z$ such 
that $K_X+B\sim _{\mathbb Q}M+\pi^*N$. 
We take a log minimal model of $(X, B)$ over $Z$ by using 
the arguments in \cite[Section 3]{birkar}. 
Then we obtain a projective surjective morphism $\pi_Y:Y\to Z$ such that 
$(Y/Z, B_Y+\sum _j E_j)$ is a log minimal 
model of $(X/Z, B)$, 
where $B_Y$ is the pushforward of $B$ on $Y$ by 
$\phi:X\dashrightarrow Y$ and 
$E_j$ is exceptional 
over $X$ and is a prime divisor 
on $Y$ for every $j$. 
Let $A$ be a sufficiently ample general Cartier divisor 
on $Z$. 
Then $(Y, B_Y+E+\pi_Y^*A)$, where 
$E=\sum _j E_j$, is a log minimal 
model of $(X, B+\pi^*A)$ 
by Lemma \ref{leng}. 
Since $\kappa (Y, K_Y+B_Y+E+\pi_Y^*A)\geq 1$, 
$K_Y+B_Y+E+\pi_Y^*A$ is 
semi-ample by Theorem \ref{thm3}. 
In particular, $$K_Y+B_Y+E=K_Y+B_Y+E
+\pi_Y^*A-\pi_Y^*A$$ is $\pi_Y$-semi-ample. 
Thus, $$
R(Y/Z, K_Y+B_Y+E)=
\bigoplus _{m\geq 0}\pi_{Y*}\mathcal O_Y(\llcorner 
m(K_Y+B_Y+E)\lrcorner)$$ is 
a finitely generated $\mathcal O_Z$-algebra. 
Therefore, $$R(X/Z, K_X+B)=\bigoplus _{m\geq 0} 
\pi_*\mathcal O_X(\llcorner m(K_X+B)\lrcorner)$$ is 
a finitely generated $\mathcal O_Z$-algebra by Lemma \ref{iso}. 
We finish the proof.  
\end{proof}

The final theorem in this section is a generalization of Fukuda's theorem 
(see Theorem \ref{thm3}). 

\begin{thm}[A special case of the log 
abundance theorem]\label{final}
Let $\pi:X\to Z$ be a proper surjective morphism from a four-dimensional 
log canonical 
pair $(X, B)$ such that 
$B$ is an effective $\mathbb Q$-divisor 
and that $K_X+B$ is $\pi$-nef. When $Z$ is a point, we 
further assume that 
$\kappa (X, K_X+B)>0$. Then 
$K_X+B$ is $\pi$-semi-ample.  
\end{thm} 
\begin{proof}
Without loss of generality, we can assume that 
$\pi$ has connected fibers. 
By Proposition \ref{prop-a} 
and the log abundance theorem in dimension 
$\leq 3$, $K_{X_\eta}+B_\eta$ is nef and abundant, where 
$\eta$ is the generic point of $Z$. By Theorem \ref{main2}, 
$\bigoplus _{m\geq 0} \pi_*\mathcal O_X(\llcorner m(K_X+B)\lrcorner)$ is a finitely generated $\mathcal O_Z$-algebra. 
Therefore, $K_X+B$ is $\pi$-semi-ample by Lemma \ref{lem-b} below. 
\end{proof}

The next lemma is well known. We leave the proof 
as an exercise for the reader. 

\begin{lem}[{cf.~\cite[Theorem 2.3.15]
{lazarsfeld}}]\label{lem-a} 
Let $\pi: X\to Z$ be a projective surjective 
morphism from 
a smooth variety $X$ to a normal 
variety $Z$ and $M$ a $\pi$-nef and 
$\pi$-big Cartier divisor on $X$. 
Then $\bigoplus _{m\geq 0} \pi_*\mathcal O_X(mM)$ is a finitely 
generated $\mathcal O_Z$-algebra if and 
only if $M$ is $\pi$-semi-ample. 
\end{lem}

By \cite[Proposition 6-1-3]{kmm}, we can reduce Lemma \ref{lem-b} 
to Lemma \ref{lem-a}. 

\begin{lem}\label{lem-b} 
Let $\pi:X\to Z$ be a proper surjective 
morphism between normal varieties 
and $M$ a $\pi$-nef and $\pi$-abundant 
Cartier divisor on $X$. 
Then $\bigoplus _{m\geq 0} \pi_*\mathcal O_X(mM)$ is a finitely generated $\mathcal O_Z$-algebra if and only if 
$M$ is $\pi$-semi-ample. 
\end{lem}

\subsection{Appendix}
In this appendix, we explicitly 
state the results in dimension $\leq 3$ 
because we can find no good references 
for the relative statements (cf.~\cite{fujita}, 
\cite{kemm}, and \cite{kemm2}). 

\begin{thm}\label{aa}
Let $\pi:X\to Z$ be a proper surjective 
morphism between normal varieties. 
Assume that $(X, B)$ is log canonical 
with $\dim X\leq 3$ and 
that $B$ is an effective 
$\mathbb Q$-divisor. 
Then $$\bigoplus _{m\geq 0}\pi_*\mathcal O_X
(\llcorner m(K_X+B)\lrcorner)$$ is a finitely 
generated $\mathcal O_Z$-algebra. 
\end{thm}

\begin{proof}
When $Z$ is a point, this theorem is well known 
(cf.~\cite{fujita}, \cite{kemm}, and \cite{kemm2}). 
So, we assume that $\dim Z\geq 1$. By 
the arguments in the proof of Theorem \ref{main}, 
we can prove that $\bigoplus _{m\geq 0} 
\pi_*\mathcal O_X(\llcorner m(K_X+B)\lrcorner)$ is 
a finitely generated $\mathcal O_Z$-algebra. 
\end{proof}

\begin{thm}\label{bb} 
Let $\pi:X\to Z$ be a proper surjective 
morphism such that $(X, B)$ is log canonical 
with $\dim X\leq 3$. 
Assume that $K_X+B$ is $\pi$-nef and 
$B$ is an effective $\mathbb Q$-divisor. 
Then $K_X+B$ is $\pi$-semi-ample. 
\end{thm}

\begin{proof}
When $Z$ is a point, this theorem is well known 
(cf.~\cite{fujita}, \cite{kemm}, and \cite{kemm2}). 
So, we 
assume that $\dim Z\geq 1$. Without loss of generality, we 
can assume that $\pi$ has connected fibers. 
It is well known that $K_{X}+B$ is $\pi$-nef and 
$\pi$-abundant by 
the log abundance theorem 
in dimension $\leq 2$. 
By Theorem \ref{aa} and Lemma \ref{lem-b}, 
$K_X+B$ is $\pi$-semi-ample. 
\end{proof}

We close this appendix with a remark. 

\begin{rem}\label{cc}
Let $\pi:X\to Z$ be a proper surjective 
morphism between normal varieties. 
Assume that $(X, B)$ is klt and 
that $B$ is an effective 
$\mathbb Q$-divisor. 
Then $$\bigoplus _{m\geq 0}\pi_*\mathcal O_X
(\llcorner m(K_X+B)\lrcorner)$$ is a finitely 
generated $\mathcal O_Z$-algebra by 
\cite{bchm}. 
Therefore, by Lemma \ref{lem-b}, 
$K_X+B$ is $\pi$-semi-ample if 
and only if 
$K_X+B$ is $\pi$-nef and 
$\pi$-abundant by Lemma \ref{lem-b}. 

Of course, we know that $K_X+B$ is $\pi$-semi-ample 
if and only if 
$K_X+B$ is $\pi$-nef and 
$\pi$-abundant without appealing 
\cite{bchm}. See, for example, 
\cite{fujino-kawamata}. 
It is known as Kawamata's 
theorem (cf.~\cite{kawamata2}). 
\end{rem}

\section{Abundance theorem for irregular varieties}\label{sec4}

In this section, we treat the abundance conjecture for 
irregular varieties. 
Let us recall the following 
minimal model conjecture. 

\begin{conj}[Minimal model conjecture]\label{conj1}
Let $X$ be a smooth projective 
variety. Assume that 
$K_X$ is pseudo-effective. 
Then there exists a normal projective  
variety $X'$ which satisfies the following 
conditions{\em{:}} 
\begin{itemize}
\item[(i)] $X'$ is birationally equivalent to 
$X$. 
\item[(ii)] $X'$ has only $\mathbb Q$-factorial 
terminal singularities. 
\item[(iii)] $K_{X'}$ is nef. 
\end{itemize}
We call $X'$ a {\em{minimal model}} 
of $X$. 
\end{conj}
In Conjecture \ref{conj1}, if $K_{X'}$ is 
semi-ample, $X'$ is 
usually called a {\em{good minimal model}} 
of $X$. 

\begin{conj}[Abundance conjecture]\label{conj11} 
Let $X$ be a projective 
variety with only canonical 
singularities. Assume that 
$K_X$ is nef. Then $K_X$ is semi-ample. 
In particular, $\kappa (X)=\kappa (X, K_X)$ is non-negative. 
\end{conj}
We know that Conjectures \ref{conj1} and 
\ref{conj11} hold in dimension $\leq 3$ 
(cf.~\cite{kmm}, \cite{fafa}, etc.). 

\begin{rem}
In Conjecture \ref{conj1}, by \cite{bchm}, 
we can replace (ii) with the following 
slightly weaker condition:~(ii$'$) $X'$ has 
at most canonical singularities. 
Similarly, we can assume that $X$ has only $\mathbb Q$-factorial 
terminal singularities in Conjecture \ref{conj11}. 
\end{rem}

\begin{rem}
Let $X$ be a smooth projective 
variety. Then $X$ is uni-ruled if and only if 
$K_X$ is not pseudo-effective by \cite{bdpp}. 
\end{rem}

The next theorem is the main theorem of this section. 

\begin{thm}[Abundance theorem for irregular 
$(n+1)$-folds]\label{thm-abun}
Assume that {\em{Conjectures \ref{conj1}}} and {\em{\ref{conj11}}} hold 
in dimension $\leq n$. 
Let $X$ be a normal complete 
$(n+1)$-fold with only canonical singularities. 
If 
$K_X$ is nef and 
$q(X)\ne 0$, then $K_X$ is semi-ample.  
\end{thm}

\begin{proof}
Let $\pi:\overline X\to X$ be a resolution and 
$\alpha:\overline X\to A=\Alb (\overline X)$ 
the Albanese mapping. 
By the assumption, we have $\dim A\geq 1$. 
Since $X$ has only rational 
singularities, 
$\beta=\alpha \circ \pi^{-1}:X\to A$ is a morphism 
(cf.~\cite[Proposition 2.3]{reid}, \cite[Lemma 2.4.1]{bs})

\begin{cla}\label{cla1}
We have $\kappa (X, K_X)=\kappa (\overline X, 
K_{\overline X})\geq 0$.
\end{cla}
\begin{proof}[Proof of {\em{Claim \ref{cla1}}}] 
Let $f:\overline X\to S$ be the Stein factorization 
of $\alpha$ and $F$ a general fiber of $f$. 
Then, by \cite[Corollary 1.2]{kawamata2.5}, we have 
$$
\kappa (\overline X, K_{\overline X})\geq 
\kappa (F, K_F)+\kappa (\overline S, K_{\overline S}), 
$$ 
where 
$\overline S$ is a resolution 
of $S$. We note that 
$\kappa (\overline S, K_{\overline S})\geq 0$ 
because $S\to \beta (X)\subset A$ is generically finite 
(see, for example, \cite[Theorem 6.10, Lemma 10.1]{ueno}). 
We also note that 
$\kappa (F, K_F)=\kappa (G, K_G)\geq 0$ since 
$\dim G\leq n$, 
$G$ has only canonical singularities, 
and $K_G$ is nef, where 
$G=\pi(F)$. 
Here, we used Conjectures \ref{conj1} and \ref{conj11} in dimension 
$\dim G=\dim F\leq n$. 
Therefore, we obtain $\kappa (\overline X, K_{\overline X})\geq 0$. 
\end{proof}
\begin{cla}\label{cla2}
If $\kappa (X, K_X)=0$, 
then $\nu(X, K_X)=0$. 
\end{cla}
\begin{proof}[Proof of {\em{Claim \ref{cla2}}}]
By Kawamata's theorem (cf.~\cite[Theorem 1]{kawamata1}), 
$\beta$ is surjective and 
$\beta_*\mathcal O_X\simeq \mathcal O_A$. 
Let $G$ be a general fiber of 
$\beta$. Then 
$\kappa (G, K_G)=0$ 
by $$0=\kappa (X, K_X)\geq \kappa (G, K_G)+\kappa (A, K_A)=
\kappa (G, K_G)$$  
as in Claim \ref{cla1} and 
$\kappa (G, K_G)\geq 0$ by Conjecture \ref{conj11} in $\dim G\leq n$ since 
$K_G$ is nef. 
We note that $X$ and $G$ have only canonical 
singularities. 
By Remark \ref{cc} and the assumption:~Conjecture \ref{conj11}  
in dimension $\leq n$, 
$K_X$ is $\beta$-semi-ample. 
Therefore, 
$\beta:X\to A$ can be written as follows: 
\begin{align*}
\beta: X\overset{f}\longrightarrow S\overset{g}\longrightarrow A, 
\end{align*}
where $K_X\sim _{\mathbb Q}f^*D$ for some $g$-ample 
$\mathbb Q$-Cartier $\mathbb Q$-divisor 
$D$ on $S$, 
$g:S\to A$ is a birational morphism, 
and $S$ is a normal variety. Since 
$\kappa (X, K_X)=0$, we 
obtain $\kappa (S, D)=0$. 
So, it is sufficient to prove that $D\sim _{\mathbb Q}0$. 
By \cite[Theorem 0.2]{ambro-comp}, 
we can write $D\sim _{\mathbb Q}K_S+\Delta_S$ such that 
$(S, \Delta_S)$ is klt. 
In particular, $\Delta_S$ is effective. 
By Lemma \ref{lem-a2} below, 
we obtain 
that $g$ is an isomorphism. 
Therefore, $D\sim _{\mathbb Q}0$ since 
$\kappa (S, D)=0$ and $S=A$ is an Abelian variety. 
\end{proof}
By Claim \ref{cla1} and Claim \ref{cla2}, 
$\nu (X, K_X)>0$ implies 
$\kappa (X, K_X)>0$. In this case, 
we obtain $\kappa (X, K_X)=\nu(X, K_X)$ by Kawamata's 
argument and the assumption:~Conjectures 
\ref{conj1} and \ref{conj11} in dimension $\leq n$ 
(see the proof of \cite[Theorem 7.3]{kawamata2}). 
Therefore, $K_X$ is semi-ample by Remark \ref{cc}. 
\end{proof}

We already used the following lemma in the proof of 
Claim \ref{cla2}. 

\begin{lem}\label{lem-a2} 
Let $g:S\to A$ be a projective 
birational 
morphism from a klt pair $(S, \Delta_S)$ 
to an Abelian variety $A$. Assume that 
$K_S+\Delta_S$ is $g$-nef and 
$\kappa (S, K_S+\Delta_S)=0$. Then 
$g$ is an isomorphism. 
\end{lem}
\begin{proof}
By replacing $S$ with its small projective $\mathbb Q$-factorialization 
(cf.~\cite{bchm}), 
we can assume that $S$ is $\mathbb Q$-factorial. 
We note that $K_S=E$, where 
$E$ is effective and $\Supp E=\Exc (g)$ since 
$A$ is an Abelian variety. 
If $B=g_*\Delta_S\ne 0$, then 
$g^*B\leq m(K_S+\Delta_S)$ for some 
$m>0$. 
In this case, 
$$1\leq \kappa (A, B)=\kappa (S, g^*B)\leq 
\kappa (S, K_S+\Delta_S)=0.$$ 
It is a contradiction. 
Therefore, $B=0$. 
This means that 
$\Delta_S$ is $g$-exceptional. 
Thus, $K_S+\Delta_S$ is effective, $g$-exceptional, 
and $\Exc (g)=\Supp (K_S+\Delta_S)$. By the assumption, $K_S+\Delta_S$ is $g$-nef. 
So, $g$ is an isomorphism by the negativity lemma. 
\end{proof}

As a special case of Theorem \ref{thm-abun}, 
we obtain the abundance theorem for irregular fourfolds. 

\begin{cor}[Abundance theorem for irregular fourfolds]\label{cor45}
Let $X$ be a normal complete fourfold with only canonical 
singularities. 
Assume that $K_X$ is nef and the irregularity $q(X)$ 
is not zero. 
Then $K_X$ is semi-ample. 
\end{cor}
\begin{proof}
It is obvious by Theorem \ref{thm-abun} because 
Conjectures \ref{conj1} and \ref{conj11} hold in dimension 
$\leq 3$ (cf.~\cite{fafa}, 
\cite{km}, etc.). 
\end{proof}

We close this section with the following theorem. 

\begin{thm}[Good minimal models of irregular fourfolds]\label{thm-good}
Let $X$ be a smooth 
projective irregular  
fourfold. 
If $X$ is not uni-ruled, then $X$ has 
a good minimal model. 
More precisely, 
there exists a normal 
projective variety $X'$ such that 
$X'$ has only $\mathbb Q$-factorial 
terminal singularities, $X'$ is birationally equivalent to 
$X$, and $K_{X'}$ is semi-ample.  
\end{thm}
\begin{proof}
We run the 
minimal model program. 
Then we obtain a minimal model 
$X'$ of $X$ since $K_X$ is pseudo-effective 
by the assumption. 
Here, we used the existence and 
the termination 
of four-dimensional terminal flips (cf.~\cite[Theorem 5-1-15]{kmm}, 
\cite{shokurov}, and 
\cite[Corollary 5.1.2]{hm}). We note 
that $q(X')=h^1(X', \mathcal O_{X'})=
q(X)\ne 0$. 
Therefore, by Theorem \ref{thm-abun}, 
we obtain that $K_{X'}$ is semi-ample. 
\end{proof}

\subsection{Appendix} 
In this appendix, we give a remark on the abundance conjecture for 
fourfolds for the reader's convenience. 

\begin{conj}[Abundance conjecture for fourfolds]\label{conj-ap} 
Let $X$ be a complete fourfold 
with only canonical singularities. 
If $K_X$ is nef, then $K_X$ is semi-ample. 
\end{conj}

This conjecture is still open. 
By Corollary \ref{cor45} and Kawamata's argument 
(cf.~\cite[Theorem 7.3]{kawamata2}), 
we can reduce Conjecture \ref{conj-ap} to the following 
two problems. 

\begin{prob}\label{prob1} 
Let $X$ be a smooth projective fourfold. If $X$ is not uni-ruled 
and $q(X)=0$, then $\kappa (X)\geq 0$. 
\end{prob}

\begin{prob}\label{prob2} 
Let $X$ be a projective fourfold with only $\mathbb Q$-factorial 
terminal singularities. 
If $K_X$ is nef, $q(X)=0$, and 
$\kappa (X, K_X)=0$, then $K_X$ is numerically trivial, equivalently, 
$K_X\sim _{\mathbb Q}0$. 
\end{prob}

We explain the reduction argument closely. 
Let $X$ be a complete fourfold with only canonical 
singularities such that 
$K_X$ is nef. 
If $q(X)\ne 0$, then $K_X$ is semi-ample by Corollary \ref{cor45}. 
So, from now on, 
we can assume that $q(X)=0$. 
By taking a resolution of $X$ and running the minimal model 
program (cf.~\cite[Theorem 5-1-15]{kmm}, 
\cite{shokurov}, and 
\cite[Corollary 5.1.2]{hm}), 
there exists a projective 
variety $X'$ such that 
$K_{X'}$ is nef and that 
$X'$ has only $\mathbb Q$-factorial terminal singularities. 
Let 
$$
X\overset{f}\longleftarrow W\overset{g}\longrightarrow X'
$$
be a common resolution. Then 
$f^*K_X=g^*K_{X'}$ by the negativity lemma. 
Therefore, we can replace $X$ with $X'$ in order to prove Conjecture 
\ref{conj-ap}. 
If we solve Problem \ref{prob1}, 
then we obtain $\kappa (X, K_X)\geq 0$ since 
$X$ has only terminal singularities. 
Furthermore, if we solve Problem \ref{prob2}, 
then we can prove that $\nu(X, K_X)>0$ implies $\kappa (X, K_X)>0$. 
By the proof of \cite[Theorem 7.3]{kawamata2}, we obtain 
$\nu (X, K_X)=\kappa (X, K_X)$ 
(cf.~Proposition \ref{prop-a}). 
Thus, $K_X$ is semi-ample (cf.~Remark \ref{cc}). 

\ifx\undefined\bysame
\newcommand{\bysame|{leavemode\hbox to3em{\hrulefill}\,}
\fi

\end{document}